\documentclass{article}
\usepackage[utf8]{inputenc}
\usepackage{appendix}
\usepackage{float}
\usepackage{enumerate}
\usepackage{todonotes}
\usepackage{comment}

\title{A note on the first Steklov eigenvalue on planar domains}
\author{
Azahara DelaTorre\footnote{Sapienza Universit\`a di Roma, Dipartimento Matematica Guido Castelnuovo,  Piazzale Aldo Moro 5, 00185 Roma Italy
Email:\, azahara.delatorrepedraza@uniroma1.it},
Gabriele Mancini
\footnote{Gabriele Mancini, Dipartimento di Matematica, Università degli Studi di Bari Aldo Moro, Via Orabona 4, 70125 Bari Italy
Email:\, gabriele.mancini@uniba.it},
Angela Pistoia\footnote{Sapienza Universit\`a di Roma\,, Dipartimento di Scienze di Base e Applicate per l'Ingegneria\,, Via Scarpa 16\,, 00161 Roma\,, Italy. Email:\, angela.pistoia@uniroma1.it}   and Luigi Provenzano\footnote{Sapienza Universit\`a di Roma\,, Dipartimento di Scienze di Base e Applicate per l'Ingegneria\,, Via Scarpa 16\,, 00161 Roma\,, Italy. Email:\, luigi.provenzano@uniroma1.it}}
\date{\today}

\usepackage{makeidx}
\usepackage{amssymb , amsmath, amsthm}
\usepackage{graphicx}
\usepackage{xcolor}
\usepackage{cancel}
\newtheorem{defi}{Definition} 
\newtheorem{thm}[defi]{Theorem}

 \newtheorem{lem}[defi]{Lemma}

 \newcommand{\eps}{\epsilon}

\newcommand{\matrice}{\begin{pmatrix}}
\newcommand{\ok}{\end{pmatrix}}

\newcommand{\R}{\mathbb R}

\begin{document}

\maketitle

\noindent
{\bf Abstract.} We consider the first positive Steklov eigenvalue on planar domains. First, we provide an example of a planar domain for which a first eigenfunction has a closed nodal line. Second, we establish a lower bound for the first positive eigenvalue on certain symmetric domains and show that this eigenvalue is simple for all ellipses. These results complement two statements contained in \cite{KS}.
\vspace{11pt}


\vspace{6pt}

\section{Introduction and statement of the main results}

Let $\Omega$ be a bounded connected open set (briefly, a bounded domain) of $\mathbb R^n$. We consider the Steklov problem
\begin{equation}\label{steklov_c}
\begin{cases}
\Delta u=0\,, & {\rm in\ }\Omega\,,\\
\partial_{\nu}u=\sigma u\,, & {\rm on\ }\partial\Omega.
\end{cases}
\end{equation}
Here and in what follows, $\nu$ denotes the unit outer normal to $\partial\Omega$. If the domain is sufficiently regular (e.g., Lipschitz is enough), then problem \eqref{steklov_c} admits an increasing sequence of non-negative eigenvalues of finite multiplicity diverging to $+\infty$:
$$
0=\sigma_0<\sigma_1\leq\sigma_2\leq\cdots\leq\sigma_k\leq\cdots\nearrow+\infty.
$$
The eigenvalue $\sigma_0=0$ has multiplicity one and the corresponding eigenfunctions are the constant functions on $\Omega$. The first positive eigenvalue is $\sigma_1$. When the domain is not smooth, \eqref{steklov_c} is understood in a weak sense. Namely, a function $u\in H^1(\Omega)$ is a Steklov eigenfunction with eigenvalue $\sigma\in\mathbb R$ if
$$
\int_{\Omega}\nabla u\cdot\nabla\phi=\sigma\int_{\partial\Omega}u\phi\,,\ \ \ \forall\phi\in H^1(\Omega).
$$
 Here $H^1(\Omega)$ denotes the usual Sobolev space of functions in $L^2(\Omega)$ with weak first order derivatives in $L^2(\Omega)$.

A sequence  $\{u_k\}_{k=0}^{\infty}$ of eigenfunctions of \eqref{steklov_c} can be chosen in such a way that their traces form an orthonormal basis of $L^2(\partial\Omega)$. We will be interested in the first positive eigenvalue $\sigma_1$. The eigenvalues admit  the well-known variational characterization
$$
\sigma_k=\min_{\substack{0\ne u\in H^1(\Omega)\\\int_{\partial\Omega}uu_j=0\,, \forall j=0,...,k-1}}\frac{\int_{\Omega}|\nabla u|^2}{\int_{\partial\Omega}u^2}.
$$
In particular

\begin{equation}\label{minmax1}
\sigma_1=\min_{\substack{0\ne u\in H^1(\Omega)\\\int_{\partial\Omega}u=0}}\frac{\int_{\Omega}|\nabla u|^2}{\int_{\partial\Omega}u^2}.
\end{equation}

In this note we mainly  focus on the planar case, i.e., $n=2$.  This work is motivated by \cite{KS}, where the authors establish a lower bound for $\sigma_1$ on certain planar symmetric regions, and a Courant-type result for Steklov eigenfunctions. Though the final results are correct, part of the statements contain some flaws.

\medskip

The statement of the Courant's Theorem for Steklov eigenfunction in \cite[\S 3]{KS} reads as follows
\begin{thm}[\cite{KS} \S 3]
The nodal lines of $u_k$ divide $\Omega$ in no more than $k$ subdomains, and no nodal line is a closed curve.
\end{thm}
The second part of the statement is not correct without adding a missing assumption. Indeed, in \cite{KS}, the main argument in the proof of the Courant's Theorem is that if a nodal line is a closed curve, it implies the existence of a subdomain of $\Omega$ where $u_k$ vanishes at the boundary, which is not possible being $u_k$ harmonic. This statement is true if $\Omega$ is simply connected. If $\Omega$ is multiply connected (i.e., its boundary has more than two connected components), nothing prevents a nodal line to be homotopic to some of the boundary components. The final result of \cite{KS}, the Courant's nodal-line theorem, holds true since what the authors really use is the fact that the closure of any nodal domain has non-empty intersection with the boundary of $\Omega$, see e.g., \cite{KKP}. 

To see in a simple way that a Steklov eigenfunction can have an interior nodal line one needs to consider the Steklov problem on a Riemannian surface. For example, consider the flat cylinder $\mathbb S^1\times(-T,T)$ with coordinates $(\theta,z)\in[0,2\pi)\times(-T,T)$. It is well-known that for $T>T^*$, $\sigma_1$ is simple and a corresponding eigenfunction is $u_1(\theta,z)=z$. Here $T^*$ is the unique positive zero of the equation $z\tanh(z)=1$. See \cite{Pol17} for the computation of the Steklov spectrum of Riemannian cylinders. Hence the nodal line is a closed curve - the circle $\mathbb S^1\times\{0\}$ - that does not meet the boundary of the cylinder. One may ask if this can happen also for planar domains. If one considers a planar annulus $A(r_0):=\{r_0^2<x^2+y^2<1\}$, $r_0\in(0,1)$, it is not difficult to show that a first eigenfunction, in polar coordinates $(r,\theta)$ of $\mathbb R^2$, is of the form $u_1(r,\theta)=(a_{r_0} r+\frac{b_ {r_0}}{r})\cos(\theta+\alpha)$, for suitable $a_{r_0},b_{r_0}\in\mathbb R$ and for any $ \alpha\in\R$ (the multiplicity of the eigenvalue is $2$). See e.g., \cite{dittmar} for the computation of the Steklov spectrum on planar annuli. Therefore, all nodal lines meet the two boundary components. However, we can produce an example of a doubly connected planar domain (i.e.,  diffeomorphic to an annulus) for which the nodal line of a first eigenfunction $u_1$ is a closed curve homotopic to the boundary components. This is the first result of this paper.

\begin{thm}\label{thm_annulus}
There exists a doubly connected, planar domain $\Omega$ for which a first Steklov eigenfunction $u_1$ has a closed nodal line homotopic to the boundary components.
\end{thm}

Next, we come to the lower bound on $\sigma_1$ established in \cite[\S4, Formula (5)]{KS}.

\begin{thm}[\cite{KS} \S4, Formula (5)]\label{thm_KS}
Let $\Omega$ be a smooth planar domain. Assume that there exists a Cartesian coordinate system $(x,y)$ and smooth functions $f:(-a,a)\to(0,+\infty)$ and $g:(-b,b)\to(0,+\infty)$ such that $\Omega=\{(x,y):-f(x)<y<f(x)\,,-a<x<a\}=\{(x,y):-g(y)<x<g(y)\,,-b<y<b\}$, for some $a,b>0$. Then
\begin{equation}\label{eq_low_KS}
\sigma_1\geq\frac{1}{\max\{\|f\sqrt{1+f'^2}\|_{\infty},\|g\sqrt{1+g'^2}\|_{\infty}\}}.
\end{equation}
\end{thm}

In particular, we are interested in the application of Theorem \ref{thm_KS} to the ellipse $\mathcal E:\frac{x^2}{a^2}+\frac{y^2}{b^2}<1$ in \cite[\S 5]{KS}. The lower bound \cite[\S5, Formula (9)]{KS} reads $\sigma_1\geq\frac{1}{\max\{a,b\}}$. If we suppose $a>b$, this just reads $\sigma_1\geq \frac{1}{a}$. But this is not possible since $\sigma_1$ (and all the eigenvalues) become arbitrarily small as $b\to 0$, for fixed $a$. 
To see this, just take $u_0=x$ in the Rayleigh quotient \eqref{minmax1} for $\sigma_1$. This is possible because $\int_{\partial\mathcal E}u_0=0$, by symmetry. Thus, from \eqref{minmax1} we get
$$
	\sigma_1\leq\frac{\int_{\mathcal E}|\nabla u_0|^2}{\int_{\partial\mathcal E}u_0^2}\leq\frac{|\mathcal E|}{\int_{\partial\mathcal E}u_0^2}=\frac{\pi ab}{\int_{\partial\mathcal E}u_0^2}.
$$
	Now, an easy computation shows that $\int_{\partial \mathcal E}u_0^2\to2\int_{-a}^au_0^2=4a^3/3$ as $b\to 0$. Therefore we obtain 
	$$
	\sigma_1\leq\frac{3\pi b}{4a^2}\to 0, \text{ as }   b\to 0.
	$$
We point out that the mistake in \cite[\S5, Formula (9)]{KS} comes from the interpretation of the lower bound \eqref{eq_low_KS} in the case of the ellipse. Indeed, in the case of the ellipse $\mathcal E$, we may assume that $a>b$ and we have
\begin{equation}\label{fg_ell}
f(x)=b\sqrt{1-\frac{x^2}{a^2}},\quad \text{ and } \quad g(y)=a\sqrt{1-\frac{y^2}{b^2}}.
\end{equation}

We observe that $\frac{d}{dx}\left[f^2(1+f'^2)\right]=\frac{2b^2(b^2-a^2)x}{a^4}<0$ for $x\in(0,a)$, hence the maximum of $f^2(1+f'^2)$, and thus of $f\sqrt{1+f'^2}$  is attained at $x=0$: $$
\|f\sqrt{1+f'^2}\|_{\infty}=f(0)\sqrt{1+f'(0)^2}=b.
$$
Analogously, $\frac{d}{dy}\left[g^2(1+g'^2)\right]=\frac{2a^2(a^2-b^2)y}{b^4}>0$ for $y\in(0,b)$, hence the supremum is given by
$$
\|g\sqrt{1+g'^2}\|_{\infty}=\lim_{y\to b^-}g(y)\sqrt{1+g'(y)^2}=\frac{a^2}{b}.
$$
Thus,  inequality \eqref{eq_low_KS} just reads $\sigma_1\geq\frac{b}{a^2}$ (instead of $\sigma_1\geq\frac{1}{\max\{a,b\}}$ as stated in \cite[\S5, Formula (9)]{KS}). The lower bound is now consistent with the upper bound $\frac{3\pi b}{4a^2}$ found above.

\medskip

Rephrasing in a suitable way the proof of  \cite[\S 4, Formula (5)]{KS}, we obtain a more general result for domains of revolution in $\mathbb R^n$, which in particular implies Theorem \ref{thm_KS}.
\begin{thm}\label{thm_correct_b}
Let $\Omega$ be a smooth bounded domain in $\mathbb R^n,$ $n\geq 2$. Assume that there exists a Cartesian coordinate system $(\bar x,x_n)$ in $\mathbb R^n$ where $\bar x=(x_1,...,x_{n-1})$, and a smooth function $f:(-a,a)\to(0,+\infty)$, $f(\pm a)=0$, such that $\Omega=\{(\bar x,x_n):-a<x_n<a\,,|\bar x|<f(x_n)\}$, for some $a>0$. Let $\sigma>0$ be a Steklov eigenvalue with an associated eigenfunction which is odd with respect to some $x_i$, $i=1,...,n-1$. Then
\begin{equation}\label{eq_ndim}
\sigma\geq\frac{1}{\|f\sqrt{1+f'^2}\|_{\infty}}.
\end{equation}
If $\Omega$ is a smooth, bounded, planar domain and there exists a Cartesian coordinate system $(x,y)$ in $\mathbb R^2$ and smooth functions $f:(-a,a)\to(0,+\infty)$ and $g:(-b,b)\to(0,+\infty)$ such that $\Omega=\{(x,y):-f(x)<y<f(x)\,,-a<x<a\}=\{(x,y):-g(y)<x<g(y)\,,-b<y<b\}$, for some $a,b>0$, then the eigenspace corresponding to $\sigma_1$ has dimension at most $2$. Moreover, in the eigenspace corresponding to $\sigma_1$ there is at most one  eigenfunction that is even in $x$ and odd in $y$, and at most one eigenfunction which is odd in $x$ and even in $y$, and the following lower bound holds:
\begin{equation}\label{eq_corr}
\sigma_1\geq\min\left\{\frac{1}{\|f\sqrt{1+f'^2}\|_{\infty}},\frac{1}{\|g\sqrt{1+g'^2}\|_{\infty}}\right\}.
\end{equation}
In particular, if there is an eigenfunction associated to $\sigma_1$ which is even in $x$ and odd in $y$, we have
\begin{equation}\label{eq_x}
	\sigma_1\geq\frac{1}{\|f\sqrt{1+f'^2}\|_{\infty}}.
	\end{equation}
Similarly, if there is an eigenfunction associated to $\sigma_1$ which is odd in $x$ and even in $y$, the bound reads
	\begin{equation}\label{eq_y}
		\sigma_1\geq\frac{1}{\|g\sqrt{1+g'^2}\|_{\infty}}.
	\end{equation}
\end{thm}

By inspecting the proof of Theorem \ref{thm_correct_b} we note that, in the case of the ellipse $\mathcal E:\frac{x^2}{a^2}+\frac{y^2}{b^2}<1$ with $a>b$ and $f,g$ given by \eqref{fg_ell},  we have the strict inequality in \eqref{eq_corr}, since $f(x)\sqrt{1{+}f'(x)^2}=\|f\sqrt{1+f'^2}\|_{\infty}=b$ only for $x=0$,  and $g(y) \sqrt{1+g'(y)^2} < \|g \sqrt{1+{g'}^2}\|_\infty  $  for any $y\in (-b,b)$. Note that \eqref{eq_corr} is equivalent to \eqref{eq_low_KS}.

\medskip

As a consequence of \eqref{eq_corr}, \eqref{eq_x} and \eqref{eq_y} we have our third result concerning the simplicity of the first Steklov eigenvalues on ellipses.
\begin{thm}
The first eigenvalue $\sigma_1$ of any ellipse (different from a disk) is simple.
\end{thm}
From Theorem \ref{thm_correct_b} we deduce that, in the case of the ellipse $\mathcal E:\frac{x^2}{a^2}+\frac{y^2}{b^2}<1$, the eigenspace corresponding to $\sigma_1$ has dimension at most $2$, and it is spanned either by an eigenfunction even in $x$ and odd in $y$, or by an eigenfunction odd in $x$ and even in $y$, or by two eigenfunctions, one of each of those types. If $a>b$, from \eqref{eq_x} we deduce that that an eigenvalue $\sigma$ associated with an eigenfunction even in $x$ and odd in $y$ must satisfy $\sigma>\frac{1}{b}$. 

Now, we invoke the  Weinstock's inequality \cite{weinstock}: $\sigma_1\leq\frac{2\pi}{|\partial\Omega|}$, holding for any simply connected planar domain. If $\Omega$ is not a disk, the inequality is strict. In the case $\Omega=\mathcal E$ with $a>b$, we have that $|\partial\mathcal E|\geq 2\pi b$, that is, the perimeter of $\mathcal E$ is strictly larger that the one of the disk of radius $b$ (which is inscribed in $\mathcal E$). This implies that $\sigma_1<\frac{1}{b}$. Therefore an eigenfunction which is even in $x$ and odd in $y$ cannot be a second eigenfunction. We deduce that the second eigenvalue is simple with eigenfunction odd in $x$ and even in $y$.

\medskip

The next sections are devoted to the proofs of Theorems \ref{thm_annulus} and \ref{thm_correct_b}. More precisely, we prove Theorem \ref{thm_annulus} in Section \ref{sec:annulus} and Theorem \ref{thm_correct_b} in Section \ref{sec:bound}.

\section{Proof of Theorem \ref{thm_annulus}}\label{sec:annulus}

\begin{figure}[htp]
\centering
\includegraphics[width=0.24\textwidth]{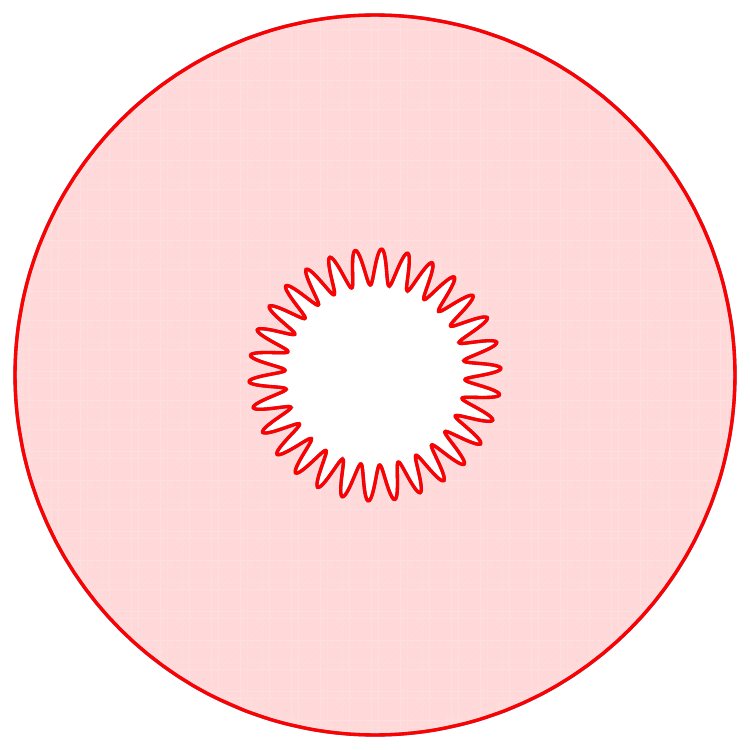}
\caption{A planar annulus with rapidly oscillating inner boundary}
\label{fig_planar}

\end{figure}

 Let $A(r_0)$ be the  annulus in the plane bounded by the concentric circles centered at the origin of radii $r_0\in(0,1)$ and $1$: $A(r_0)=\{x\in\mathbb R^2:r_0<|x|<1\}$. As we have mentioned in the introduction, for all $r_0\in(0,1)$, a first Steklov eigenfunction changes sign on each of the boundary components, hence its nodal lines touch the boundary (in four points). 
In order to find a domain having a first eigenfunction with a nodal line  entirely contained in the interior of the domain, we will perturb the inner boundary of $A(r_0)$ ($|x|=r_0$)  in a rapidly oscillating way.
 The new boundary will be contained in a small neighborhood of $|x|=r_0$ (say, of size $\epsilon$ with $\epsilon>0$ small), but its length will be comparable to $2\pi c r_0$, with $c>1$:  the length of the perturbed inner boundary will be larger than the length of the original inner boundary. We will then choose the factor $c>0$ in a suitable way. See Figure \ref{fig_planar}.

To be more precise in defining the deformed annulus, we need to recall an  auxiliary result.

\begin{lem}\label{lemma_bucur}
Let $\Omega$ be a smooth domain in $\mathbb R^n$ and let $w\in L^{\infty}(\partial\Omega)$, $w\geq 1$ be a density on $\partial\Omega$. Then there exists a sequence of domains $\Omega_{\eps}$ that converges uniformly to $\Omega$ as $\epsilon\to 0$ and such that $\mathcal H|_{\partial\Omega_{\eps}}\rightharpoonup w\mathcal H|_{\partial\Omega}$
weakly* in the sense of measures, where $\mathcal H$ is the $n-1$ dimensional Hausdorff measure. In particular,
$$
|\partial\Omega_{\eps}|\to\int_{\partial\Omega}w$$
and
$$
\sigma_k(\Omega_{\eps})\to\sigma_k^w(\Omega),
$$
where $\sigma_k^w(\Omega)$ are the eigenvalues of the following weighted Steklov problem on $\Omega$
\begin{equation}\label{weighted_O}
\begin{cases}
\Delta u=0\,, & {\rm in\ }\Omega\,,\\
\partial_{\nu}u=\sigma w u\,, & {\rm on\ }\partial\Omega 
\end{cases}
\end{equation}
and $\sigma_k(\Omega_{\epsilon})$ are the Steklov eigenvalues  on $\Omega_{\eps}$. Moreover, if $\sigma_k^w(\Omega)$ is a simple eigenvalue of \eqref{weighted_O}, then there exists $\epsilon_0>0$ such that for all $\epsilon\in(0,\epsilon_0)$, $\sigma_k(\Omega_{\epsilon})$ is simple, and if $u_k^w$ is an eigenfunction associated with $\sigma_k^w(\Omega)$ with $\int_{\partial\Omega}w(u_k^w)^2=1$,  there exists a subsequence $u_k^{\epsilon_j}$,  $\epsilon_j\to 0$ as $j\to\infty$, of eigenfunctions associated with $\sigma_k(\Omega_{\epsilon_j})$, with $\int_{\partial\Omega_{\epsilon_j}}(u_k^{\eps_j})^2=1$ such that
$$
\lim_{j\to\infty}\|u_k^{\epsilon_j}-u_k^w\|_{H^1(\Omega_{\delta}\cap\Omega\cap\Omega_{\eps})}=0
$$
for all $\delta>0$, where $\Omega_{\delta}:=\{x\in\Omega:d(x,\partial\Omega)>\delta\}$. Here $d(\cdot,\partial\Omega)$ denotes the distance to the boundary of $\Omega$.
\end{lem}
One can find a detailed proof e.g., in \cite{bucur,FL_Steklov}. Note that for any $\delta>0$, we can choose $\epsilon'(\delta)$ small enough, so that for all $\eps<\eps'(\delta)$ we have $\Omega_{\delta}\cap\Omega\cap\Omega_{\eps}=\Omega_{\delta}$, since the convergence of the domains is uniform. In \cite{FL_Steklov} the statement on the convergence of eigenfunctions is given also for the case of multiple eigenvalues. However, the statement is more involved and requires the introduction of the concept of convergence of generalized eigenfunctions. Since, as we shall see, we will treat simple eigenvalues,  we will need the result as it is stated here above. Following \cite{FL_Steklov}, if one wants a precise limit weight $w\geq 1$, for example, a constant weight $c$, in the case of a two-dimensional domain, one can produce a quite explicit perturbation, for example, considering as boundary of $\partial\Omega_{\epsilon}$ the graph of a periodic, rapidly oscillating function on the circles defining $\partial\Omega$.

\medskip

We consider now the following weighted Steklov problem in $A(r_0)$:

\begin{equation}\label{weigh_AB}
\begin{cases}
\Delta u=0 & {\rm in\ } A(r_0)\\
\partial_{\nu}u=\sigma w u  & {\rm on\ }\partial A(r_0),
\end{cases}
\end{equation}
where $w\equiv 1$ at $|x|=1$ and $w\equiv c\geq 1$ at $|x|=r_0$ is a positive weight on $\partial A(r_0)$. We denote the eigenvalues by
$$
\sigma_k^w(A(r_0)).
$$
We choose $c=\frac{1}{r_0}$, and solve problem \eqref{weigh_AB}. 
 It is  standard to verify, by separation of variables,  that any eigenfunction can be written  in polar coordinates $(r,\theta)$ as $u(r,\theta)=(a_lr^l+b_lr^{-l}){ \cos(l \theta+\alpha)}$, for suitable $a_l,b_l$ and any $\alpha\in\mathbb R$ (when $l\geq 1$), or it is of the form $u(r,\theta)=a_0+b_0\log(r)$ (when $l=0$).
The boundary condition involving the eigenvalue on $\partial A(r_0)=\{x: |x|=1\}\cup\{x: |x|=r_0\}$ is determined  by the relations:
$$
\begin{cases}
\partial_ru(1,\theta)=\sigma u(1,\theta),\\
-\partial_r u(r_0,\theta)=\frac{\sigma}{r_0}u(r_0,\theta).
\end{cases}
$$
This corresponds to a system of two linear equations in two unknowns $a_l,b_l$ which has a non-trivial solution if and only if the determinant of the associated matrix vanished. Such matrix is 
$$
	\begin{pmatrix}
		l-\sigma & -l-\sigma\\
		r_0^{-1+l}(-\sigma-l) &  r_0^{-1-l}(l-\sigma)
	\end{pmatrix}
	$$
when $l\geq 1$, and
$$
\begin{pmatrix}
	-\sigma & 1\\
-	\frac{\sigma}{r_0} & - \frac{1+\sigma\log(r_0)}{r_0}
\end{pmatrix}
$$
when $l=0$.
Imposing the vanishing of the determinant we obtain the two solutions 
$$
\sigma_-(l)=l\frac{1-r_0^l}{1+r_0^l}\,,\ \ \ \sigma_+(l)=l\frac{1+r_0^l}{1-r_0^l}
$$
for $l\geq 1$, and two solutions when $l=0$:
$$
\sigma_-(0)=0\,,\ \ \ \sigma_+(0)=-\frac{2}{\log(r_0)}.
$$
All these numbers exhaust the spectrum of \eqref{weigh_AB}. Let $R^*\in(0,1)$ be the unique positive root of $-\frac{2}{\log(r)}=\frac{1-r}{1+r}$. It is a calculus exercise to see that $\sigma_1^w(A(r_0))=\sigma_-(1)$ for $r_0<R^*$ and it is double, while $\sigma_1^w(A(r_0))=\sigma_+(0)$ for $r_0>R^*$ and it is simple. For $r_0=R^*$ the eigenvalue is triple. From now on we will choose $r_0>R^*$. Hence the first eigenvalue is

$$
\sigma_1^w(A(r_0))=-\frac{2}{\log(r_0)},
$$
it is simple, and a corresponding eigenfunction is
$$
u_1^w(r,\theta)=\kappa\left(1-\frac{2\log(r)}{\log(r_0)}\right)
$$ 
which is radial. Here $\kappa=\frac{1}{2\sqrt{\pi}}$, so that $\int_{\partial\Omega}w(u_1^w)^2=1$.  Hence the nodal line is defined by the equation $r=\sqrt{r_0}>r_0$. This is not surprising. In fact, the eigenvalues of \eqref{weigh_AB} are just the eigenvalues of the usual Steklov problem on the flat cylinder $\mathbb S^1\times(-T,T)$ with $T=-\frac{1}{2}\log(r_0)$, and the equation defining $R^*$ corresponds, through the conformal map from the cylinder to the annulus $A(r_0)$ defined by $(\theta,z)\mapsto e^{-(z+L+i\theta)}$, to the equation defining $T^*$, that is $z\tanh(z)=1$. Recall that $T^*$ is the critical value of the cylinder height for which we have a triple eigenvalue, and a change of angular momentum for the first eigenfunction. What we are truly doing, is approximating a Steklov problem on a cylinder, written as a weighted Steklov problem for an annulus through a conformal map, with classical Steklov problems defined on planar domains with oscillating boundaries.

\medskip

We apply now Lemma \ref{lemma_bucur} to $A(r_0)$ with the weight $w\equiv\frac{1}{r_0}$ on the inner boundary and $w\equiv 1$ on the outer boundary. We deduce that there exists a sequence of domains $\Omega_{\eps}(r_0)$ such that
$$
	\lim_{\eps\to 0}|\partial\Omega_{\eps}(r_0)|=\int_{\partial A(r_0)}w=4\pi
	$$
and
$$
\lim_{\eps\to 0}\sigma_1(\Omega_{\eps}(r_0))=\sigma_1^w(A(r_0))=-\frac{2}{\log(r_0)}.
$$
Moreover, up to taking a subsequence, if $u^{\eps}_1$ denotes the first Steklov eigenfunction on $\Omega_{\eps}$ satisfying $\int_{\partial\Omega_{\eps}}(u_1^{\eps})^2=1$ we have that, for any $\delta>0$,
$$
\lim_{\eps\to 0}\|u_1^{\eps}-u_1^w\|_{H^1(A(r_0)_{\delta})}=0
$$
where $A(r_0)_{\delta}=\{x\in\mathbb R^2:r_0+\delta<|x|<1-\delta\}$. Then, by standard elliptic regularity estimates (see e.g., \cite[Thm. 8.10]{GT}) we find that 
$$
\lim_{\eps\to 0}\|u_1^{\eps}-u_1^w\|_{C^1(A(r_0)_{\delta})}=0.
$$
Now, at the points where $u_1^w=0$ we have $\nabla u_1^w\ne 0$. We conclude by Thom's isotopy Theorem (see e.g., \cite[\S 20.2]{AbRo}) that for $\varepsilon>0$ small enough, $u_1^{\eps}$ has a nodal line homotopic to the nodal line of $u_1^w$, which is the circle $r=\sqrt{r_0}$.

\section{Proof of Theorem \ref{thm_correct_b}}\label{sec:bound}

\begin{figure}[htp]
\centering
\includegraphics[width=0.9\textwidth]{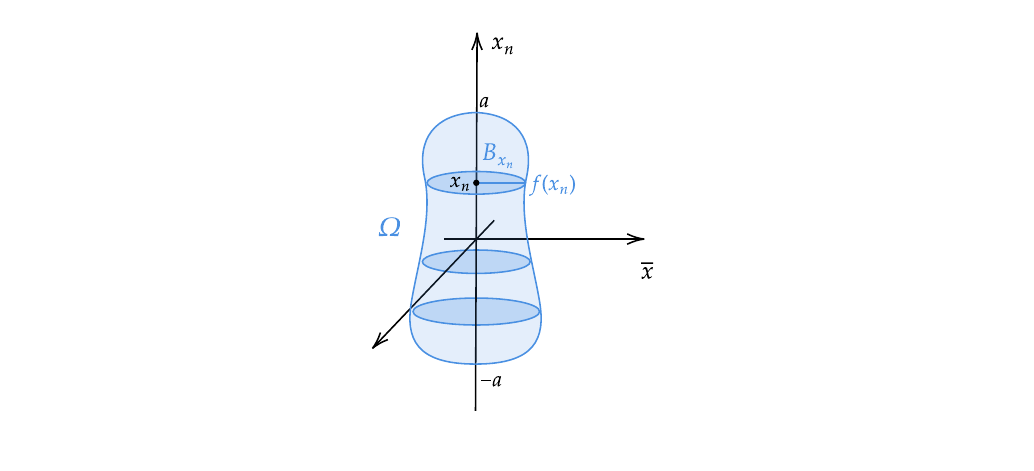}
\caption{A domain of revolution foliated by balls.}
\label{foliation}
\end{figure}

Let $\Omega$ be a smooth domain of revolution as in the statement of Theorem \ref{thm_correct_b}, namely
$$
\Omega=\{(\bar x,x_n)\in\mathbb R^n:-a<x_n<a\,, |\bar x|<f(x_n)\},
$$
for a smooth function $f>0$ in $(-a,a)$ and  with $f(\pm a)=0$. 
Let $z\in \R$ and let $B_z:=\Omega\cap\{x_n=z\}$ be the $n-1$ dimensional ball of radius $f(z)$ obtained by intersecting $\Omega$ with the hyperplane $x_n=z$. We have that $\Omega$ is foliated by $n-1$-dimensional balls: $\Omega=\bigcup_{x_n\in(-a,a)}B_{x_n}$, see Figure \ref{foliation}.

Let now $u$ be a Steklov eigenfunction associated with some eigenvalue $\sigma$, which is odd with respect to some $x_i$, $i=1,...,n-1$. { For $x_n\in (-a,a),$} consider $u|_{B_{x_n}}$. It follows  {  by symmetry} that $\int_{\partial B_{x_n}}u=0$; hence,
$$
\frac{\int_{B_{x_n}}|\nabla_{\bar x}u|^2}{\int_{\partial B_{x_n}}u^2}\geq\sigma_1(B_{x_n})=\frac{1}{f(x_n)}.
$$
Here $\nabla_{\bar x}$ denotes the gradient with respect to $\bar x$ and $\sigma_1(B_{x_n})$ is the first (positive) Steklov eigenvalue of the $n-1$ dimensional ball $B_{x_n}$ of radius $f(x_n)$, which is $\frac{1}{f(x_n)}$. We have then
$$
		\int_{B_{x_n}}|\nabla u|^2\geq\int_{B_{x_n}}|\nabla_{\bar x}u|^2\geq\frac{1}{f({x_n})}\int_{\partial B_{x_n}}u^2\frac{\sqrt{1+f'({x_n})^2}}{\sqrt{1+f'({x_n})^2}}
		\geq\frac{1}{\|f\sqrt{1+f'^2}\|_{\infty}}\int_{\partial B_{x_n}}u^2\sqrt{1+f'({x_n})^2}.
$$
Integrating both sides in ${x_n}\in(-a,a)$ we get
$$
\int_{\Omega}|\nabla u|^2\geq\frac{1}{\|f\sqrt{1+f'^2}\|_{\infty}}\int_{\partial\Omega}u^2.
$$
This proves \eqref{eq_ndim}.

\medskip

Let now $\Omega$ be a smooth, bounded, planar domain as in the second part of the statement. Let $u$ be a Steklov eigenfunction with eigenvalue $\sigma$ which is odd with respect to $y$. Then, the proof above with $x_n=x$, $\bar x=y$ gives that
\begin{equation}\label{low1}
\sigma\geq\frac{1}{\|f\sqrt{1+f'^2}\|_{\infty}}.
\end{equation}
Note that the basic inequality used in the first part of the proof, namely,
$$
\frac{\int_{B_{x_n}}|\nabla_{\bar x}u|^2}{\int_{\partial B_{x_n}}u^2}\geq\sigma_1(B_{x_n})=\frac{1}{f(x_n)},
$$
in this case just reads
\begin{equation}\label{fund0}
\frac{\int_{-f(x)}^{f(x)}(\partial_yu(x,y))^2dy}{u^2(x,-f(x))+u^2(x,f(x))}\geq\frac{1}{f(x)}
\end{equation}
which can be deduced also from the Fundamental Theorem of Calculus applied to an odd function on a symmetric interval. More precisely, if $u(x,y)$ is odd in $y$ for all $x$,
$$
u(x,f(x))=u(x,f(x))-u(x,0)=\int_0^{f(x)}\partial_yu
$$
which implies
\begin{equation}\label{fund1}
u^2(x,f(x))=\left(\int_0^{f(x)}\partial_yu\right)^2\leq f(x)\int_0^{f(x)}(\partial_yu)^2
\end{equation}
and analogously, by parity,
\begin{equation}\label{fund2}
u^2(x,-f(x))\leq f(x)\int_{-f(x)}^0(\partial_yu)^2.
\end{equation}
Summing \eqref{fund1} and \eqref{fund2} gives \eqref{fund0}.
Inequality \eqref{fund0} is just telling us that the first positive Steklov eigenvalue of the segment $(-f(x),f(x))$ is larger or equal to $\frac{1}{f(x)}$. It is actually equal since we know that for a segment $(-L,L)$ the Steklov problem has exactly two eigenvalues: $0$, $1/L$.

\medskip

Analogously, if $\sigma$ is a Steklov eigenvalue with eigenfunction which is odd with respect to $x$, we get, by the same argument above with $x_n=y$, $\bar x=x$, that
\begin{equation}\label{low2}
\sigma\geq\frac{1}{\|g\sqrt{1+g'^2}\|_{\infty}}.
\end{equation}

Finally, we note that, being $\Omega$ symmetric with respect to both $x$ and $y$ by construction, any eigenfunction $u$ can be written as a combination of symmetric eigenfunctions. More precisely, $u=u_{ee}+u_{oe}+u_{eo}+u_{oo}$. Here $u_{ee}$ denotes a function which is even in $x$ and $y$, $u_{oe}$ a function which is odd in $x$ and even in $y$,   $u_{eo}$ a function which is even in $x$ and odd in $y$, and  $u_{oo}$ a function which is odd in $x$ and  $y$. If we take $u=u_1$, that is, a first eigenfunction, we see that necessarily $u_{ee}=u_{oo}=0$, since by Courant's Theorem, any first eigenfunction has exactly two nodal domains, and this is not possible if the function is even or odd in both $x$ and $y$. Moreover if we have two linearly independent eigenfunctions both of the type $u_{oe}$ (or both of the type $u_{eo}$), since they must be orthogonal between them, and orthogonal to the constant, we get a contradiction with having two nodal sets. Hence the first eigenspace has dimension at most $2$, and a first eigenfunction must be odd with respect to one of the axes. Hence inequalities \eqref{low1} and \eqref{low2} hold with $\sigma=\sigma_1$, and in particular, they imply
$$
\sigma_1\geq\min\left\{\frac{1}{\|f\sqrt{1+f'^2}\|_{\infty}},\frac{1}{\|g\sqrt{1+g'^2}\|_{\infty}}\right\}.
$$





\section*{Acknowledgments}
The first, second and third authors are members of the Gruppo Nazionale per l'Analisi  Matematica, la Probabilit\`{a} e le loro Applicazioni (GNAMPA) of the Istituto Nazionale di Alta Matematica (INdAM) {and were supported by the 
INdAM-GNAMPA Project 2024, codice CUP E53C23001670001, and the INdAM-GNAMPA Project 2025, codice CUP E5324001950001}. The fourth author acknowledges the support of the INdAM GNSAGA group. The first author  acknowledges financial support from the Spanish Ministry of Science and Innovation (MICINN), through the IMAG-Maria de Maeztu Excellence Grant CEX2020-001105-M/AEI/10.13039/501100011033. She is also supported by the FEDER-MINECO Grants PID2021- 122122NB-I00 and PID2020-113596GB-I00; RED2022-134784-T, funded by MCIN/AEI/10.13039/501100011033 and by J. Andalucia (FQM-116); Fondi Ateneo – Sapienza Università di Roma; PRIN (Prot. 20227HX33Z). 
The second author thanks the PNRR MUR project CN00000013 HUB - National Centre for HPC, Big Data and Quantum Computing (CUP H93C22000450007) and acknowledges financial support from the PRIN PNRR  P2022YFAJH project {\lq\lq Linear and Nonlinear PDEs: New directions and applications''} (CUP H53D23008950001).
The fourth author aknowledges financial support from the project ``Perturbation problems and asymptotics for elliptic differential equations: variational and potential theoretic methods'' funded by the European Union – Next Generation EU and by MUR-PRIN-2022SENJZ3.

\section*{Statements and Declarations}
{\bf Conflict of Interest}
The authors declare that they have no conflict of interest.

{\bf Data Availability}
This article does not contain any data that requires sharing.

\bibliography{bibliography}{}
\bibliographystyle{abbrv}
\end{document}